\documentclass[12pt]{amsart}

\usepackage{amssymb}
\usepackage{enumerate}

\newtheorem{thm}{Theorem}[section]

\newtheorem{prop}[thm]{Proposition}

\newtheorem{cor}[thm]{Corollary}

\newenvironment{rk}{\refstepcounter{thm}\noindent%

{\bf Remark \arabic{section}.\arabic{thm}} \ }{

\smallskip

}

\newenvironment{pf}[1][]{\noindent {\it Proof #1} : }{\hbox{~}\qed \smallskip}

\def\Ndb{\mathbb N}

\def\Rdb{\mathbb R}

\def\Zdb{\mathbb Z}

\newcommand{\vfi}{\varphi}
\def\la{\lambda}
\def \norm{| \! |}

\title{Embeddings of proper metric spaces into Banach spaces}
\author{Florent BAUDIER$^{\dag}$}

\begin{document}

\maketitle

\begin{abstract}
We show that there exists a strong uniform embedding from any proper metric space into any Banach space without cotype. Then we prove a result concerning the Lipschitz embedding of locally finite subsets of $\mathcal{L}_{p}$-spaces. We use this locally finite result to construct a coarse bi-Lipschitz embedding for proper subsets of any $\mathcal{L}_p$-space into any Banach space $X$ containing the $\ell_p^n$'s. Finally using an argument of G. Schechtman we prove that for general proper metric spaces and for Banach spaces without cotype a converse statement holds. 
\end{abstract}

\makeatletter

\renewcommand{\@makefntext}[1]{#1}

\makeatother

\footnotetext{ Laboratoire de Math\'ematiques, UMR 6623

 Universit\'e de Franche-Comt\'e,

 25030 Besan\c con, cedex - France \\ \indent $^{\dag}$florent.baudier@univ-fcomte.fr\\
 \indent \subjclass{2000 {\it Mathematics Subject Classification}.}{ (46B20)
 (51F99)}}

\section{Introduction and Notation}
We fix some notation and give some motivation concerning the embedding of metric spaces.

Let $(M,d)$ and $(N,\delta)$ be two metric spaces and $f:M\to N$ be a map. For
$t>0$ define :
$$\rho_f(t)={\rm Inf}\{\delta(f(x),f(y)),\ \ d(x,y)\geq t\}$$
and
$$\omega_f(t)={\rm Sup}\{\delta(f(x),f(y)),\ \ d(x,y)\leq t\}.$$

We say that $f$ is a {\it coarse embedding} if $\omega_f(t)$ is finite for all
$t>0$ and $\lim_{t\to \infty}\rho_f(t)=\infty$.

Suppose now that $f$ is injective. We say that $f$ is a {\it uniform embedding} if $f$ and $f^{-1}$ are uniformly continuous (i.e $\omega_f(t)$ and $\omega_{f^{-1}}(t)$ tend to $0$ when $t$ tends to $0$). 
According to \cite{K} we refer to $f$ as a {\it strong uniform embedding} if it is simultaneously a uniform embedding and a coarse embedding. 

Following \cite{MN},
we also define the {\it distortion} of $f$ to be
$$ {\rm dist}(f):= \|f\|_{Lip}\|f^{-1}\|_{Lip}=\sup_{x\neq y \in
M}\frac{\delta(f(x),f(y))}{d(x,y)}.\sup_{x\neq y \in
M}\frac{d(x,y)}{\delta(f(x),f(y))}.$$ If the distortion of $f$ is finite, we
say that $f$ is a {\it Lipschitz embedding} and that $M$ Lispchitz embeds into
$N$. The terminology {\it metric embedding} is used too.

In the sequel our terminology for classical notions in Banach space theory follows \cite{LT1} and \cite{AK}.
Let $X$ and $Y$ be two Banach spaces. If $X$ and $Y$ are linearly isomorphic,
the {\it Banach-Mazur distance} between $X$ and $Y$, denoted by $d_{BM}(X,Y)$,
is the infimum of $\|T\|\,\|T^{-1}\|$, over all linear isomorphisms $T$ from
$X$ onto $Y$.

For $p\in [1,\infty]$ and $n\in \Ndb$, $\ell_p^n$ denotes the space $\Rdb^n$
equipped with the $\ell_p$ norm. We say that a Banach space $X$ uniformly
contains the $\ell_p^n$'s if there is a constant $C \geq 1$ such that for every
integer $n$, $X$ admits an $n$-dimensional subspace $Y$ so that
$d_{BM}(\ell_p^n,Y)\leq C$.

The problem of embedding of metric spaces occurs very often in different fields of mathematics. Let us focus on the following two.
The first one deals with Lipschitz classification of Banach spaces and is related to the celebrated result of I. Aharoni \cite{A} which asserts that every separable metric space Lipschitz embeds into $c_0$ (the Banach space of all real sequences converging
to $0$, equipped with the supremum norm). Is the converse statement of I. Aharoni's result true? Namely, if a Banach space which is Lipschitz universal for the class of all separable metric spaces does it necessarily contain a isomorphic copy of $c_0$?
This question is equivalent to the following problem :

Problem I : If $c_0$ Lipschitz embeds into a Banach space $X$ does it imply that $c_0$ is isomorphic to a closed subspace of $X$ ?

This question is difficult but a positive result for the isometric case can be found in \cite{GK} where G. Godefroy and N. J. Kalton
showed, using Lipschitz-free Banach spaces, that if a separable Banach space is isometric to a subset of another Banach space $Y$, then it is actually linearly isometric to a subspace of $Y$. There exist counterexamples in the nonseparable case. Very recently, in \cite{DL} Y. Dutrieux and G. Lancien prove that a Banach space containing isometrically every compact metric space must contain a subspace linearly isometric to $C([0,1])$ (the space of all continuous functions on $[0,1]$ equipped with the supremum norm). In the same vein, in \cite{K2} N. J. Kalton proved that if $c_0$ coarsely or uniformly embeds into a Banach space $X$ then one of its iterated duals has to be nonseparable (in particular, $X$ cannot be reflexive). On the other hand, N. J. Kalton proved in \cite{K} that $c_0$ can be strongly uniformly embedded into a Banach space with the Schur property. Thus $c_0$ can be uniformly and coarsely embedded into a Banach space which does not contain a linear copy of $c_0$.

The generalization of Index Theory for non-compact Riemmanian manifolds gives birth to the second problem.

Problem II: Which metric spaces can be coarsely embedded into a super-reflexive or
reflexive Banach space?

The strong result of Kasparov and Yu in \cite{KY} gives a precise idea of this connection.

They showed that if $M$ is a discrete metric space with bounded geometry and if $M$ is coarsely embeddable into a uniformly convex Banach space, then the coarse geometric Novikov conjecture holds for $M$.

The definition of a metric space with bounded geometry is the following :

A metric space $M$ is {\it locally finite} if any ball of $M$ with finite
radius is finite. If moreover, there is a function $C:(0,+\infty)\to \Ndb$ such
that any ball of radius $r$ contains at most $C(r)$ points, we say that $M$ has
a {\it bounded geometry}.

Problem II received an upsurge of interest in the light of Kasparov and Yu's result.
Some significant partial answers have been given in \cite{BG} and \cite{BL} where it is proved that any locally finite metric space Lipschitz embeds into any Banach space without cotype. We recall that a metric space is {\it stable} if for any pair of sequences $(x_n)_{n=1}^{\infty}$, $(y_n)_{n=1}^{\infty}$ such that both iterated limits exist we have $$\lim_{m\to\infty}\lim_{n\to\infty}d(x_m,y_n)=\lim_{n\to\infty}\lim_{m\to\infty}d(x_m,y_n).$$ In \cite{K2} it is proved that a stable metric space can always be coarsely or uniformly embedded into a reflexive Banach space. In the next section we improve this result for proper metric spaces, which are automatically stable.

\section{Strong uniform embedding of proper metric spaces}
A metric space is proper if all its closed balls are compact.
In this section we will prove the following theorem :
\begin{thm}\label{strong}
 Let $(M,d)$ a proper metric space and $X$ a Banach space without cotype, then $M$ strongly uniformly embeds into $X$. 
\end{thm}

\begin{pf}
Let $X$ be a Banach space uniformly containing the $\ell_{\infty}^n$'s (or equivalently without any nontrivial cotype).

Fix $t_0\in M$, a proper metric space, and denote for $n\in\Zdb$,\\ $B_n=B(t_0,2^{n+1})$ the closed ball of $M$ of radius $2^{n+1}$ centered at the point $t_0$.

Let $M_n^k$ a maximal $2^{-k+n+3}$-net of $B_n$, with $k\ge 1$. 

We define $$\begin{array}{rcl}

     \vfi_n^k\ :\ B_n  & \rightarrow & \ell_\infty(M_n^k):=X_{n,k}\\
       &   &  \\
     t & \mapsto & \big(d(t,s)-\vert s\vert\big)_{s\in M_n^k},\ {\rm where}\ \vert s\vert=d(s,t_0). \\
    \end{array}$$

Let $\Phi$: $I=\Zdb\times\Ndb\to\Zdb^+$ a bijection.
Following \cite{BL} we can build inductively finite dimensional subspaces
$(F_j)_{j=0}^\infty$ of $X$ and $(T_j)_{j=0}^\infty$ so that for every $j\geq
0$, $T_j$ is a linear isomorphism from $X_{\Phi^{-1}(j)}$ onto $F_j$ satisfying
$$\forall u\in X_{\Phi^{-1}(j)}\ \ \ \ \frac{1}{2}\|u\|\leq \|T_ju\|\leq \|u\|$$
and also such that $(F_j)_{j=0}^\infty$ is a Schauder finite dimensional
decomposition of its closed linear span $Z$. Let $P_j$ be the
projection from $Z$ onto $F_0\oplus...\oplus F_j$ with kernel $\overline{\rm
Span}\,(\bigcup_{i=j+1}^\infty F_i)$, we may assume after renorming, that $\|P_j\|\leq
1$. Finally we denote $\Pi_0=P_0$ and $\Pi_j=P_j-P_{j-1}$ for $j\geq 1$. We have
that $\|\Pi_j\|\leq 2$. 

Define now $$\begin{array}{rcl}

     f_n^k\ :\ B_n  & \rightarrow & F_{\Phi(n,k)}\\
       &   &  \\
     t & \mapsto & T_{\Phi(n,k)}(\vfi_n^k(t)) \\
    \end{array}$$

and,  $$\begin{array}{rcl}

     f_n\ :\ B_n  & \rightarrow &\displaystyle \sum_{k\ge 1} \Pi_{\Phi(n,k)}(Z)\\
       &   &  \\
     t & \mapsto &\displaystyle \sum_{k\ge 1}\frac{f_n^k(t)}{(n-k)^2+1} \\
    \end{array}$$
It is obvious that $f_n$ is $C$-Lipschitz with $C=\displaystyle\sum_{n\in\Zdb}\frac{1}{n^2+1}$.\\
Finally the embedding is $$\begin{array}{rcl}

    f\ :\ M & \rightarrow & Z=\overline{\rm{Span}}(\bigcup_{i=0}^{\infty}F_{i})\subset X \\
            &   &  \\
          t & \mapsto & \la_t f_{n}(t)+(1-\la_t)f_{n+1}(t)\ ,\ {\rm if}\ 2^{n}\le\vert t\vert\le 2^{n+1} \\
    \end{array}$$
where, $$\la_t=\frac{2^{n+1}-\vert t\vert}{2^{n}},\ n\in\Zdb$$\\

\medskip We start by showing that $f$ is Lipschitz. Let $a,b \in M$ and assume,
as we may, that $|a|\leq |b|$.

\medskip\noindent I) If $|a|\leq \frac{1}{2}|b|$, then
$$\|f(a)-f(b)\|\leq C(|a|+|b|)\leq \frac{3C}{2}|b|\leq 3C(|b|-|a|)\leq
3Cd(a,b).$$

\medskip\noindent II) If $\frac{1}{2}|b|<|a|\leq |b|$, we have two different
cases to consider.

\smallskip 1) $2^n\leq |a|\leq |b|<2^{n+1}$, for some $n$. Then,
let
$$\la_a=\frac{2^{n+1}-|a|}{2^{n}}\ \ {\rm and}\ \ \la_b=\frac{2^{n+1}-|b|}{2^{n}}.$$
We have that
$$|\la_a-\la_b|=\frac{|b|-|a|}{2^{n}}\leq \frac{d(a,b)}{2^{n}},$$
so

$$\begin{array}{rcl}
\norm f(a)-f(b)\norm &  =  & \norm\la_a f_n(a)-\la_b f_n(b)+(1-\la_a)f_{n+1}(a)\\
                     &     & -(1-\la_b)f_{n+1}(b)\norm\\
                     &     &                        \\
                     & \le & \la_a\norm f_n(a)-f_n(b)\norm +(1-\la_a)\norm f_{n+1}(a)-f_{n+1}(b)\norm\\
                     &     & +2C|\la_a-\la_b||b|\\
                     &     &                         \\
                     & \le & Cd(a,b)+2^{n+2}C|\la_a-\la_b|\\
                     &     &             \\
                     & \le & 5C\,d(a,b).
\end{array}$$

\smallskip 2) $2^n\leq |a|<2^{n+1}\leq |b|<2^{n+2}$, for some $n$. Then,
let
$$\la_a=\frac{2^{n+1}-|a|}{2^{n}}\ \ {\rm and}\ \ \la_b=\frac{2^{n+2}-|b|}{2^{n+1}}.$$
We have that,
$$\la_a \leq \frac{d(a,b)}{2^n},\ \ {\rm so}\ \ \la_a|a|\leq 2\,d(a,b).$$
Similarily,
$$1-\la_b=\frac{|b|-2^{n+1}}{2^{n+1}}\leq \frac{d(a,b)}{2^{n+1}}\ \ {\rm
and}\ \ (1-\la_b)|b|\leq 2\,d(a,b).$$

$$\begin{array}{rcl}
\|f(a)-f(b)\| &  =  & \|\la_a f_n(a)+(1-\la_a)f_{n+1}(a)-\la_b f_{n+1}(b)\\
              &     & -(1-\la_b)f_{n+2}(b)\|\\
              &     &                     \\
              & \le & \la_a(\|f_n(a)\|+\|f_{n+1}(a)\|)+(1-\la_b)(\|f_{n+1}(b)\|\\
              &     & +\|f_{n+2}(b)\|)+\|f_{n+1}(a)-f_{n+1}(b)\|\\
              &     &                     \\
              & \le & 2C\la_a|a|+2C(1-\la_b)|b|+2Cd(a,b)\\
              &     &                     \\
              & \le & 9C\,d(a,b).
\end{array}$$

\smallskip\noindent We have shown that $f$ is $9C$-Lipschitz.

\medskip We shall now estimate $f$ from below. We consider
$a,b\in M$ and assume again that $|a|\leq |b|$. We need
to study three different cases. In our discussion, whenever $|a|$ (respectively
$|b|$) will belong to $[2^{m},2^{m+1})$, for some integer $m$, we shall denote
$$\la_a=\frac{2^{m+1}-|a|}{2^{m}}\ \ ({\rm respectively}\ \ \la_b=\frac{2^{m+1}-|b|}{2^{m}}).$$

In the sequel we denote $R_{(n,k)}=T_{\Phi_{(n,k)}}^{-1}\circ\Pi_{\Phi_{(n,k)}}$ for the sake of convenience.

\smallskip 1) $2^{n}\leq |a|\leq |b|<2^{n+1}$, for some $n$.

\vskip 0.5cm
There exists $l\ge 1$ and $s_{b}\in M_n^{l+3}$ such that $2^{-l+n+2}\le d(a,b)<2^{-l+n+3}$ and $d(s_{b},b)<2^{-l+n}$. We denote $k=l+3$ in the sequel.

Thus, $R_{(n,k)}(f(a)-f(b))=\frac{1}{(n-l-3)^2+1}(\la_a\varphi_n^k(a)-\la_b\varphi_n^k(b))$
and so,
$$\begin{array}{rcl}
 \la_a\varphi_n^k(a)-\la_b\varphi_n^k(b) & = & (\la_a(d(a,s)-|s|)-\la_b(d(b,s)-|s|))_{s\in M_n^k}\\
 (\la_a\varphi_n^k(a)-\la_b\varphi_n^k(b))_{|_{s_{b}}} & = & \la_a d(a,s_{b})+(\la_b-\la_a)|s_{b}|-\la_bd(b,s_{b})\\
\end{array}$$

On the other hand, $$\begin{array}{rcl}
R_{(n+1,k)}(f(a)-f(b)) & = & \frac{1}{(n-l-2)^2+1}((1-\la_a)\varphi_{n+1}^k(a)\\
              &   &     -(1-\la_b)\varphi_{n+1}^k(b))\\
 \end{array}$$
and,
$$\begin{array}{rcl}
 (1-\la_a)\varphi_{n+1}^k(a)-(1-\la_b)\varphi_{n+1}^k(b) & = & ((1-\la_a)(d(a,s)-|s|)\\
                                                         &   & -(1-\la_b)(d(b,s)-|s|))_{s\in M_{n+1}^k}\\
 ((1-\la_a)\varphi_{n+1}^k(a)-(1-\la_b)\varphi_{n+1}^k(b))_{|_{s_{b}}} & = & (1-\la_a)d(a,s_{b})+(\la_a-\la_b)|s_{b}|\\
                                                         &   & -(1-\la_b)d(b,s_{b})\\
\end{array}$$

We remark that $$n-l-3\le n-l-2\le n-l+2\le\log(d(a,b))<n-l+3$$ thus if $0\le n-l-3$ then $$(n-l-3)^2+1\le (n-l+2)^2+1\le (\log d(a,b))^2+1$$
otherwise $$\log(d(a,b))-7\le n-l-3\le n-l-2< 0$$ and we have $$(n-l-2)^2+1\le(n-l-3)^2+1\le \left(\log\frac{d(a,b)}{2^7}\right)^2+1.$$
Let $\alpha(t)=(\log (t))^2+1$.
Now if we sum $R_{(n+1,k)}(f(a)-f(b))$ and $R_{(n,k)}(f(a)-f(b))$ we get 
$$\begin{array}{rcl}
   8\cdot{\rm Max}\left\{\alpha(d(a,b)),\alpha\left(\frac{d(a,b)}{2^7}\right)\right\}\norm f(a)-f(b)\norm & \ge & d(a,s_{b})-d(b,s_{b})\\
                                                         & \ge & d(a,b)-2d(b,s_{b})\\
							 & \ge & d(a,b)-2^{-l+n+1}\\
							 & \ge & \frac{d(a,b)}{2}
  \end{array}$$

\smallskip 2) $2^n\leq |a|<2^{n+1}\leq |b|<2^{n+2}$, for some $n$.

There exists $l\ge 1$ and $s_{a}\in M_{n+1}^{k}$ such that $2^{-l+n+3}\le d(a,b)<2^{-l+n+4}$ and $d(s_{a},a)<2^{-l+n}$.

We use again the projections and we evaluate at some specific points.

$$\begin{array}{rcl}
   R_{(n,k)}(f(a)-f(b))_{|_{t_0}} & = & \frac{1}{(n-k)^2+1}(\la_a d(t,s)-\la_a|s|)_{|_{s=t_0}}\\
 	                                                             & = & \frac{\la_a}{(n-k)^2+1}|a|\\
  \end{array}$$

$$\begin{array}{rcl}
   R_{(n+2,k)}(f(a)-f(b))_{|_{t_0}} & = & \frac{1}{(n-k+2)^2+1}(-(1-\la_b) d(b,s)\\
                                                                         &   & -(1-\la_b)|s|)_{|_{s=t_0}}\\
 & & \\
 	                                                                 & = & -\frac{(1-\la_b)}{(n-k+2)^2+1}|b|\\
  \end{array}$$

$$\begin{array}{rcl}
   R_{(n+1,k)}(f(a)-f(b))_{|_{s_a}} & = & \frac{1}{(n-k+1)^2+1}((1-\la_a) d(a,s)-(1-\la_a)|s|\\
                                                                         &   & -\la_b d(b,s)+\la_b|s|)_{|_{s=s_a}}\\
 & &\\
 	                                                                 & = &\frac{1}{(n-k+1)^2+1}((1-\la_a)d(a,s_a)-\la_bd(b,s_a)\\

                                                                         &   &-(1-\la_a)|s_a|+\la_b|s_a|\\
  \end{array}$$

Thus,

$$\begin{array}{rcl}
12{\rm Max}\left\{\alpha(d(a,b)),\alpha\left(\frac{d(a,b)}{2^7}\right)\right\}\norm f(a)-f(b)\norm  & \ge & \la_a|a|+(1-\la_b)|b|\\
                                                     &     & -(1-\la_a)d(a,s_a)+\la_bd(b,s_a)\\
                                                     &     & +(1-\la_a)|s_a|-\la_b|s_a|\\
 & & \\
		                                     & \ge & \la_bd(a,b)+\la_a|a|+(1-\la_b)|b|\\
                                                     &     & +(1-\la_a)|a|-\la_b|a|\\
                                                     &     & -\la_bd(s_a,a)-(1-\la_a)d(s_a,a)\\
                                                     &     & -(1-\la_a)d(a,s_a)-\la_bd(s_a,a)\\
 & &\\
		                                     & \ge & d(a,b)-2(1-\la_a+\la_b)d(a,s_a)\\
 & &\\
		                                     & \ge & d(a,b)-4d(s_a,a)\\
 & &\\
		                                     & \ge & \frac{d(a,b)}{2}\\
\end{array}$$

\smallskip 3) $2^n\leq|a|<2^{n+1}<2^p\leq |b|<2^{p+1}$ for some integers $n$
and $p$.

There exists $l\ge 1$ such that $2^{-l+p+2}\le d(a,b)<2^{-l+p+3}$

$$\begin{array}{rcl}
   R_{(p,l)}(f(a)-f(b))_{|_{t_0}} & = & \frac{1}{(p-l)^2+1}(-\la_b d(b,s)+\la_b|s|)_{|_{s=t_0}}\\
 	& = & -\frac{\la_b}{(p-l)^2+1}|b|\\
  \end{array}$$

$$\begin{array}{rcl}
   R_{(p+1,l)}(f(a)-f(b))_{|_{t_0}} & = & \frac{1}{(p-l+1)^2+1}(-(1-\la_b) d(b,s)\\
 & & +(1-\la_b)|s|)_{|_{s=t_0}}\\
 & &\\
 	& = & -\frac{(1-\la_b)}{(p-l+1)^2+1}|b|\\
  \end{array}$$

$$\begin{array}{rcl}
8\cdot{\rm Max}\left\{\alpha(d(a,b)),\alpha\left(\frac{d(a,b)}{2^7}\right)\right\}\norm f(a)-f(b)\norm  & \ge & \la_b|b|+(1-\la_b)|b|\\
 & &\\
		& \ge & |b| \\
 & &\\	
	        & \ge & \frac{2}{3}(|a|+|b|)\\
 & &\\
		& \ge & \frac{2}{3}d(a,b)\\
\end{array}$$

\medskip\noindent All possible cases are settled and we have shown that $f$ is a strong uniform embedding which satisfies the following estimates:
$$\gamma(d(a,b))\le\norm f(a)-f(b)\norm\le 9Cd(a,b),$$
where $$\gamma(t)=\frac{t}{24\cdot{\rm Max}\left\{\alpha(t),\alpha\left(\frac{t}{2^7}\right)\right\}}$$
has the following behavior $$\begin{array}{rcl}
       \gamma(t) & = & \frac{t}{(\log(t))^2+1}\ {\rm if}\ t>2^7\\
 & & \\
                 & \ge & D\ {\rm a\ positive\ constant\ if}\ 1\le t\le 2^7\\
 & & \\
                 & =  &  \frac{t}{\left(\log\left(\frac{t}{2^7}\right)\right)^2+1}\ {\rm if}\ t<1.\\
      \end{array}$$

\end{pf}

\section{Lispchitz embedding for locally finite subsets of $\mathcal{L}_{p}$-spaces}

Let $1\le p\le\infty$ and $\la\ge 1$. A Banach space $X$ is a $\mathcal{L}_{p,\la}$-space if for every finite-dimensional subspace $F$ of $X$ there exists a finite-dimensional subspace $G\subset X$ containing $F$ such that $d_{BM}(G,\ell_p^m)\le\la$ (where $m$ is the dimension of $G$). $X$ is a $\mathcal{L}_{p}$-space if it is a $\mathcal{L}_{p,\la}$-space for some $\la$.
We refer to \cite{BLi} Appendix F for results on $\mathcal{L}_{p}$-spaces.
Several results on embeddings of locally finite sets were obtained by M. I. Ostrovskii in \cite{O}, the author and G. Lancien in \cite{BL}.

\begin{thm}\label{baulan}[Baudier-Lancien]
 Let $X$ be a Banach space without cotype and let $(M,d)$ be a locally finite metric space. Then there exists a Lipschitz embedding of $M$ into $X$.
\end{thm}

\begin{thm}\label{ovst}[Ostrovskii]
 Let $(M,d)$ be a locally finite subset of a Hilbert space. Then $M$ is Lipschitz embeddable into an arbitrary infinite dimensional Banach space.
\end{thm}

We prove a similar result in the broad context of $\mathcal{L}_{p}$-spaces which extends theorems \ref{baulan} and \ref{ovst}.
\begin{prop}\label{script}
Let $1\le p\le\infty$ and $\la\ge 1$. Let $(Y,\norm\cdot\norm)$ a $\mathcal{L}_{p,\la}$-space and $M$ a locally finite subset of $Y$, then $M$ Lipschitz embeds into any Banach space $X$ which contains uniformly the $\ell_p^n$'s.
Moreover the distortion is universal.
\end{prop}

\begin{pf}
 Let $B_n:=\{t\in M;\norm t\norm\le 2^{n+1}\}$ for $n\in\Ndb$.
We may assume that $B(t_0,1)=\{t_0\}$ and that $t_0=0$.
$F_n=\overline{{\rm Span} B_n}$ is a finite-dimensional subspace of $Y$, hence there exists a finite-dimensional subspace $G_n$ containing $F_n$ and an isomorphism $R_n$ from $G_n$ onto $\ell_p^{dim(G_n)})$ such that ${\rm dist}(R_n)\le\la$, with $\Vert R_n\Vert\le 1$ and $\Vert R_n^{-1}\Vert\le\la$. Fix $\delta>0$.
Like in the proof of theorem \ref{strong} we build inductively a Schauder finite-dimensional decomposition
$(Z_n)_{n=0}^\infty$ of a subspace of $X$ and $(T_n)_{n=0}^\infty$ so that for every $n\geq
0$, $T_n$ is a linear isomorphism from $\ell_p^{dim(G_n)}$ onto $Z_n$ satisfying
$$\forall u\in \ell_p^{dim(G_n)}\ \ \ \ \frac{1}{1+\delta}\norm u\norm\leq \norm T_nu\norm\leq \norm u\norm.$$
We keep the same notations for the different projections.

We define $$\begin{array}{rcl}
             f_n : G_n\supset B_n & \longrightarrow & Z\subset X \\
                              t & \longmapsto & T_n\circ R_n(t) \\
            \end{array}$$

And then we construct $f:M\longrightarrow Z\subset X$ as follows:

(i) $f(0)=0$.

(ii) For $n\geq 0$ and $2^n\leq \norm t\norm<2^{n+1}$:
$$f(t)=\la_t
f_n(t)+(1-\la_t)f_{n+1}(t),\ \ {\rm where}\ \
\la_t=\frac{2^{n+1}-\norm t\norm}{2^{n}}.$$

In the sequel we will always suppose that $\norm a\norm\le\norm b\norm$. We have to consider three cases.

\begin{enumerate}
 \item $2^n\le\norm a\norm\le\norm b\norm<2^{n+1}$

We have 
\begin{equation}\label{eq1}
 \norm \Pi_n(f(a)-f(b))\norm=\norm\la_a f_n(a)-\la_b f_n(b)\norm
\end{equation}
and, \begin{equation}\label{eq2}
      \norm \Pi_{n+1}(f(a)-f(b))\norm =\norm (1-\la_a) f_{n+1}(a)-(1-\la_b) f_{n+1}(b)\norm
     \end{equation}

By the triangle inequality (\ref{eq1}) is between \begin{equation}
                                                 \la_a\norm f_n(a)-f_n(b)\norm \buildrel {+}\over{-} (\la_a-\la_b)\norm f_n(b)\norm
                                                \end{equation}

And (\ref{eq2}) between \begin{equation}
                         (1-\la_a)\norm f_{n+1}(a)-f_{n+1}(b)\norm \buildrel {+}\over{-} (\la_a-\la_b)\norm f_{n+1}(b)\norm
                        \end{equation}

So we get
$$\frac{1}{4}\left(\frac{\norm a-b\norm}{\la(1+\delta)}-2(\la_a-\la_b)\norm b\norm\right)\le\norm f(a)-f(b)\norm\le
\norm a-b\norm+2(\la_a-\la_b)\norm b\norm$$
\noindent But $(\la_a-\la_b)\norm b\norm\le\frac{\norm b\norm-\norm a\norm}{2^n}2^{n+1}$, hence if $\norm b\norm-\norm a\norm\le\frac{\norm a -b\norm}{5(1+\delta)\la}$ we have,

$$\frac{\norm a-b\norm}{20(1+\delta)\la}\le\norm f(a)-f(b)\norm\le5\norm a-b\norm$$

If $\norm b\norm-\norm a\norm\ge\frac{\norm a -b\norm}{5(1+\delta)\la}$, using the linearity we get,

$$\begin{array}{rcl}
4\norm f(a)-f(b)\norm & \ge & \frac{\norm \la_a a-\la_b b\norm}{\la(1+\delta)}+\frac{\norm (1-\la_a) a-(1-\la_b) b\norm}{\la(1+\delta)}\\
                   &     &        \\
                     & \ge &\frac{1}{\la(1+\delta)}(\la_b \norm b\norm-\la_a\norm a\norm+(1-\la_b)\norm b\norm-(1-\la_a)\norm a\norm)\\
                   &     &        \\
                     & \ge & \frac{\norm a-b\norm}{5(1+\delta)^2\la^2}\\

\end{array}$$

\item $2^n\le\norm a\norm<2^{n+1}\le\norm b\norm<2^{n+2}$

The quantity $\norm f(a)-f(b)\norm$ is between the sum and the half-max of the numbers
\begin{equation}\label{eq3}
 \norm \Pi_n(f(a)-f(b))\norm=\norm \la_a f_n(a)\norm,
\end{equation}
\begin{equation}\label{eq4}
       \norm \Pi_{n+1}(f(a)-f(b))\norm=\norm (1-\la_a) f_{n+1}(a)-\la_b f_{n+1}(b)\norm\ {\rm and}
     \end{equation}
\begin{equation}\label{eq5}
 \norm \Pi_{n+2}(f(a)-f(b))\norm=\norm (1-\la_b)f_{n+2}(b)\norm.
\end{equation}

But, $$(\ref{eq3})\le\frac{2^{n+1}-\norm a\norm}{2^n}\norm a\norm\le 2(2^{n+1}-\norm a\norm)\le 2(\norm b\norm-\norm a\norm)\le 2\norm b-a\norm.$$

$$(\ref{eq5})\le\frac{\norm b\norm-2^{n+1}}{2^{n+1}}\norm b\norm\le 2\norm b-a\norm.$$

$$\begin{array}{rcl}
(\ref{eq4}) &  =  & \norm f_{n+1}(a)-f_{n+1}(b)+(1-\la_b)f_{n+1}(b)-\la_a f_{n+1}(a)\norm\\
                   &     &        \\
            & \le & \norm a-b\norm +2\norm a-b\norm +2\norm a-b\norm\\
                   &     &        \\
            & \le & 5\norm a-b\norm.\\
\end{array}$$
On the other hand,
$$\begin{array}{rcl}
(\ref{eq3}) & \ge & \frac{\la_a\norm a\norm}{\la(1+\delta)}\\
                   &     &        \\
            & \ge & \frac{1}{\la(1+\delta)}\frac{2^{n+1}-\norm a\norm}{2^n}\norm a\norm\\
                   &     &        \\
            & \ge & \frac{2^{n+1}-\norm a\norm}{\la(1+\delta)}\\
\end{array}$$
And, $(\ref{eq5})\ge\frac{(1-\la_b)\norm b\norm}{\la(1+\delta)}\ge\frac{\norm b\norm-2^{n+1}}{\la(1+\delta)}$,
hence $$2\norm f(a)-f(b)\norm\ge{\rm Max}\left\{\frac{2^{n+1}-\norm a\norm}{\la(1+\delta)},\frac{\norm b\norm-2^{n+1}}{\la(1+\delta)}\right\}.$$
If ${\rm Max}\{2^{n+1}-\norm a\norm,\norm b\norm-2^{n+1}\}\ge\frac{\norm b-a\norm}{5\la(1+\delta)}$, then $$\norm f(a)-f(b)\norm\ge\frac{\norm b-a\norm}{10\la^2(1+\delta)^2}.$$

Otherwise, $$\begin{array}{rcl}
       (\ref{eq4}) & \ge & \norm f_{n+1}(b)-f_{n+1}(a)\norm-(1-\la_b)\norm b\norm-\la_a\norm a\norm\\
                   &     &        \\
                   & \ge & \frac{\norm b-a\norm}{\la(1+\delta)}-2(\norm b\norm-2^{n+1})-2(2^{n+1}-\norm a\norm)\\
                   &     &        \\
                   & \ge & \frac{\norm b-a\norm}{\la(1+\delta)}-\frac{4}{5}\frac{\norm b-a\norm}{\la(1+\delta)}\\
                   &     &        \\
                   & \ge & \frac{\norm b-a\norm}{5\la(1+\delta)}
      \end{array}$$
Finally, $$\frac{\norm b-a\norm}{10\la^2(1+\delta)^2}\le\norm f(a)-f(b)\norm\le9\norm a-b\norm.$$

\item $2^n\le\norm a\norm<2^{n+1}<2^p\le\norm b\norm<2^{p+1}$

We have \\$\norm \Pi_n(f(a)-f(b))\norm=\norm \la_a f_n(a)\norm$,\\ $\norm \Pi_{n+1}(f(a)-f(b))\norm=\norm (1-\la_a) f_{n+1}(a)\norm$,\\ $\norm \Pi_p(f(a)-f(b))\norm=\norm \la_b f_p(b)\norm$\\ and $\norm \Pi_{p+1}(f(a)-f(b))\norm=\norm (1-\la_b) f_{p+1}(b)\norm$.\\ Hence, $$\norm f(a)-f(b)\norm\le\norm a\norm+\norm b\norm$$ and using the projections $\Pi_p$ and $\Pi_{p+1}$, $$\norm f(a)-f(b)\norm\ge\frac{1}{4}\left(\frac{\la_b\norm b\norm}{\la(1+\delta)}+\frac{(1-\la_b)\norm b\norm}{\la(1+\delta)}\right)=\frac{\norm b\norm}{4\la(1+\delta)}.$$

But we have the following inequalities,$$\frac{\norm b\norm}{2}\le\norm b\norm-\norm a\norm\le\norm b-a\norm\le\norm b\norm+\norm a\norm\le2\norm b\norm.$$
And thus, $$\frac{\norm b-a\norm}{8\la(1+\delta)}\le\norm f(a)-f(b)\norm\le 4\norm a-b\norm.$$

\end{enumerate}

Finally we get the following global estimates: $$\frac{\norm b-a\norm}{20\la^2(1+\delta)^2}\le\norm f(a)-f(b)\norm\le 9\norm a-b\norm.$$

\end{pf}

\newpage

\begin{rk}
 We can reconstruct the proof of theorem \ref{baulan} combining proposition \ref{script} and the fact that $C([0,1])$ is a $\mathcal{L}_{\infty}$-space and isometrically universal for all separable metric spaces.

Together for $p=2$ our statement is theorem \ref{ovst}.
\end{rk}

\section{Coarse bi-Lispchitz embedding of proper metric spaces}

$f$ is a coarse bi-Lipschitz embedding if there exists two non negative constants $C_d$ and $C_a$ such that for all $x,y\in M$, $$\frac{1}{C_d}d(x,y)-C_a\le\delta(f(x),f(y))\le C_d d(x,y)+C_a.$$\\
In the sequel we will use the notation $M \buildrel{(C_d,C_a)}\over {\hookrightarrow} N$. A coarse bi-Lipschitz embedding need not be injective neither continuous. The constant $C_d$ will be called the {\it dilation coarse bi-Lipschitz constant}, and $C_a$ the {\it additive coarse bi-Lipschitz constant}.

We prefer to use the term coarse bi-Lipschitz embedding instead of quasi-isometric embedding which is used for instance in \cite{N}, because in Nonlinear Functionnal Analysis the term quasi-isometric embedding refers to a Lipschitz embedding $f$ such that ${\rm dist}(f)\le1+\epsilon$, for all $\epsilon>0$.
We remark that according to N. J. Kalton's terminology from \cite{K3}, a coarse bi-Lipschitz embedding $f$ is in particular a {\it coarse Lipschitz} map, which means that ${\rm limsup}_{t\to\infty}\frac{\omega_f(t)}{t}<\infty$.

\begin{cor}\label{coro}
 Let $1\le p\le\infty$ and $\la\ge 1$. Let $(Y,\norm\cdot\norm)$ a $\mathcal{L}_{p,\la}$-space and $M$ a proper subset of $Y$, then $M$ admits a coarse bi-Lipschitz embedding into any Banach space $X$ which contains uniformly the $\ell_p^n$'s. Moreover the additive and dilation coarse bi-Lipschitz constants are universal, they depend only on $\la$.
\end{cor}

\begin{pf}
 $B_n:=\{t\in M;\ \norm t\norm\le2^{n+1}\}$ is compact. Let $R_n$ an $\frac{\epsilon}{2}$-net of $B_n$.
The cardinal of $R_n$ is finite. $R:=\bigcup_nR_n$ is an $\frac{\epsilon}{2}$-net of $M$, and is locally finite.
Define  $$\begin{array}{rcl}

    \beta\ :\ M & \rightarrow & R \\
             t & \mapsto & \beta(t),\ {\rm such\ that\ } \norm t-\beta(t)\norm<\frac{\epsilon}{2} \\
    \end{array}$$
It is straigthforward that $\norm a-b\norm-\epsilon\le\norm\beta(a)-\beta(b)\norm\le\norm a-b\norm+\epsilon$.
If we embed $R$ using theorem \ref{script} and compose with $\beta$ then we build a coarse bi-Lipschitz embedding with universal constants $C_d$ and $C_a$.

\end{pf}

Actually in the case $p=\infty$, G. Schechtman indicated to us an argument to prove that the converse statement holds. We thank him for allowing us to present it there.
\begin{prop}
 The following assertions are equivalent :
\begin{enumerate}[i)]
 \item $X$ has no nontrivial cotype
 \item There exists a universal constant $C_d>0$ such that for every proper metric space $(M,d)$ and for all $\epsilon>0$, $M \buildrel{(C_d,\epsilon)}\over{\hookrightarrow} X$
 \item There exist two universal positive constants $C_d$ and $C_a$ such that for every proper metric space $(M,d)$, $M \buildrel{(C_d,C_a)}\over{\hookrightarrow} X$
\end{enumerate}

\end{prop}

\begin{pf}
 i) implies ii) is Corollary \ref{coro} and ii) implies iii) is trivial.

Let us prove iii) implies i) using G. Schechtman's argument.
 
Fix $n$ and $k$ be two positive integers, and let $N_k=\frac{1}{k}\Zdb^n\cap kB_{\ell_{\infty}^n}$.
$N_k$ is finite, hence there exists an embedding $\Theta_k$ : $N_k\to X$ such that $\Theta_k(0)=0$ and $N_k\buildrel{(C_d,C_a)}\over{\hookrightarrow} X$. Let $\Phi_k(x)=\frac{\Theta_k(kx)}{k}$, then $$\frac{1}{C_d}\norm x-y\norm_{\infty}-\frac{C_a}{k}\le\norm \Phi_k(x)-\Phi_k(y)\norm\le C_d \norm x-y\norm_{\infty}+\frac{C_a}{k}.$$
Define $\Psi_k$ from $B_{\ell_{\infty}^n}$ to $N_k$ such that for all $x\in B_{\ell_{\infty}^n}$, we have
 
$\norm x-\Psi_k(x)\norm_{\infty}\le\frac{1}{k}$.
Finally define $$\begin{array}{rrcl}
                \Phi : & B_{\ell_{\infty}^n} & \to     & X_{\mathcal{U}}\\
                       & x                   & \mapsto & \overline{(\Phi_k\circ\Psi_k(x))_{k\ge 1}}^{\mathcal{U}},\\
               \end{array}$$

where $X_{\mathcal{U}}$ is the ultra-product of $X$ according to a nontrivial ultrafilter $\mathcal{U}$ of $\Ndb$.
$\Phi$ is a Lipschitz embedding with distortion $C_d$ the dilation coarse bi-Lipschitz constant. Using a delicate $w*$-differentability argument due to S. Heinrich and P. Mankiewicz \cite{HM}, we can prove that $\ell_{\infty}^n$ linearly embeds into $X_{\mathcal{U}}^{**}$. And the fact that the double-dual and the ultra-product of a Banach space $Y$ are finitely representable in $Y$ allows us to conclude.

\end{pf}

\begin{rk}
 Using the argument of G. Schechtman one also can prove that the converse of theorem $2.1$ in \cite{BL} holds. Namely, a Banach space $X$ uniformly contains the $\ell_{\infty}^n$'s if and only if every locally finite metric space Lipschitz embeds into $X$.
\end{rk}


\newpage

\begin{thebibliography}{WW}

\bibitem{A} I. Aharoni, Every separable metric space is Lipschitz equivalent to
a subset of $c_0^+$, {\it Israel J. Math.}, {\bf 19} (1974), 284-291.
\\
\bibitem{AK} F. Albiac and N. J. Kalton, Topics in Banach space theory, Graduate Texts in Mathematics, vol. 233, Springer, New York, 20006.
\\
\bibitem{B} F. Baudier, Metrical characterization of super-reflexivity and linear type of Banach spaces, {\it Archiv Math.} {\bf 89} (2007), 419-429.
\\
\bibitem{Bo} J. Bourgain, The metrical interpretation of super-reflexivity in Banach spaces, {\it Israel J. Math.}  {\bf 56} (1986),
221-230.
\\
\bibitem{BL} F. Baudier and G. Lancien, Embeddings of locally finite metric spaces into Banach spaces, {\it Proc. Amer. Math. Soc.}, {\bf 136} (2008),1029-1033.
\\
\bibitem{BLi} Y. Benyamini and J. Lindenstrauss, Geometric non linear functional analysis. Vol.
1, American Mathematical Society Colloqium Publications, vol. 48, American
Mathematical Society, Providence, RI, 2000.
\\
\bibitem{BG} N. Brown and E. Guentner, Uniform embeddings of bounded geometry
spaces into reflexive Banach spaces, {\it Proc. Amer. Math. Soc.}, {\bf 133(7)}
2005, 2045-2050.
\\
\bibitem{D} A. Dvoretzky, Some results on convex bodies and Banach spaces, {\it Proc. Internat. Sympos. Linear Spaces} (Jerusalem, 1960)
123--160.
\\
\bibitem{DL} Y. Dutrieux and G. Lancien, Isometric embeddings of compact spaces into Banach spaces, {\it J. Functional Analysis}, to appear.
\\
\bibitem{GK} G. Godefroy and N. J. Kalton, Lipschitz free Banach spaces, {\it Studia Math.}, {\bf 159} (2003), 121-141.
\\
\bibitem{HM} S. Heinrich and P. Mankiewicz, Applications of ultrapowers to the uniform and Lipschitz classification of Banach spaces, {\it Studia Math.}, {\bf 73} (1982), 225-251.
\\
\bibitem{K} N. J. Kalton, Spaces of Lipschitz and Hölder functions and their applications, {\it Collect. Math.}, {\bf 55} (2004), 171-217.
\\
\bibitem{K2} N.J. Kalton, Coarse and uniform embeddings into reflexive spaces, {\it Q. J. Math.}, {\bf 58} (2007) no. 3, 393-414.
\\
\bibitem{K3} N. J. Kalton, The nonlinear geometry of Banach spaces, preprint.
\\
\bibitem{KY} G. Kasparov and G. Yu, The coarse Novikov conjecture and uniform
convexity, {\it Advances Math.}, to appear.
\\
\bibitem{LT1} J. Lindenstrauss and L. Tzafriri, Classical Banach spaces I,
Springer Berlin 1977.
\\
\bibitem{MP} B. Maurey and G. Pisier, S\'{e}ries de variables al\'{e}atoires vectorielles
ind\'{e}pendantes et propri\'{e}t\'{e}s g\'{e}om\'{e}triques des espaces de Banach, {\it Studia
Math.}, {\bf 58(1)} 1976, 45-90.
\\
\bibitem{MN} M. Mendel and A. Naor, Metric cotype, {\it Ann. Math.}, to appear.
\\
\bibitem{N} P. W. Novak, Remarks on quasi-isometric non-embeddability into uniformly convex Banach spaces,  arXiv:math/0506178v3.
\\
\bibitem{O} M. I. Ostrovskii, Coarse embeddability into Banach spaces, arXiv:0802.3666v1.
\\
\bibitem{R} M. Ribe, On uniformly homeomorphic normed spaces, {\it Ark. Mat.}, {\bf 14} (1976), 237-244.
\\
\end{thebibliography}
\end{document}